\documentclass{article}      

\usepackage{epsfig}

\title{Double q-Analytic q-Hermite Binomial Formula\\
and q-Traveling Waves}

\author{Sengul Nalci and Oktay K. Pashaev  \\Department of Mathematics, Izmir Institute of Technology \\ Urla-Izmir, 35430, Turkey}

\begin{document}
\newcommand{\be}{\begin{equation}}
\newcommand{\ee}{\end{equation}}
\newcommand{\bea}{\begin{eqnarray}}
\newcommand{\eea}{\end{eqnarray}}
\newcommand{\disp}{\displaystyle}
\newcommand{\la}{\langle}
\newcommand{\ra}{\rangle}

\newtheorem{thm}{Theorem}[subsection]
\newtheorem{cor}[thm]{Corollary}
\newtheorem{lem}[thm]{Lemma}
\newtheorem{prop}[thm]{Proposition}
\newtheorem{defn}[thm]{Definition}
\newtheorem{rem}[thm]{Remark}
\newtheorem{prf}[thm]{Proof}
\newtheorem{iden}[thm]{Identity}

\maketitle


\begin{abstract}

Motivated by derivation of the Dirac type $\delta$-function for quantum states in Fock-Bargmann representation, we find $q$-binomial expansion in terms of $q$-Hermite polynomials, analytic in two complex arguments. Based on this representation, we introduce a new class of complex functions of two complex arguments, which we call the double $q$-analytic functions.
The real version of these functions describe the $q$-analogue of traveling waves, which is not preserving the shape during evolution as the usual traveling wave. For corresponding $q$-wave equation we solve IVP in the $q$-D'Alembert form.
\end{abstract}

\section{Introduction}

In complex analysis, a complex function $f(z)$ of one complex variable $z$ is analytic in some domain $D$ if in $D$ it satisfies \begin{equation} \frac{\partial}{\partial \bar{z}} f(z)= \frac{1}{2}\left(\frac{\partial}{\partial x}+i \frac{\partial}{\partial y}\right)f(z)=0. \label{eqn0}\end{equation}

If we consider a complex function $f(z,w)$ of two complex variables $z,w$, analytic in both variables \begin{equation}\frac{\partial}{\partial \bar{z}} f(z,w)= \frac{\partial}{\partial \bar{w}} f(z,w)=0, \label{analytic2}\end{equation}
it is called a double analytic if \begin{equation} \frac{1}{2} \left( \frac{\partial}{\partial z}+ i\frac{\partial}{\partial w} \right) f(z,w)=0. \label{2analytic}\end{equation}

As an example, $f(z,w)=(z+iw)^2$ is double analytic, while $(z-iw)^2$ is double anti-analytic.

In general, double analytic function can be written as power series in complex binomials $$f(z,w)=\sum_{n=0}^\infty a_n (z+iw)^n.$$
For such double analytic binomials here we derive the following Hermite binomial formula
\begin{equation}
(z+ i w )^n=\frac{1}{2^n} \sum^n_{k=0} { n \choose k }i^k H_{n-k}(z) H_k(w) \label{eqn1} .
\end{equation}

Derivation of this formula is motivated by description of the Dirac type $\delta$- function for quantum states in the Fock-Bargmann representation.
As is well known, states of a quantum system
in the Fock-Bargmann representation  are described by complex analytic function $f(z)$ and visa versa \cite{Perelomov}.
In this representation, due to formula
\begin{equation}
\int d \mu(z) \,e^{\xi \bar{z}} f(z)=f(\xi),
\end{equation}
where measure $d \mu(z)=dz d\bar{z} e^{-z \bar{z}},$
the exponential function plays the role of the Dirac type  $\delta$ function \cite{Floratos}.
Proof of this formula is based on following identity
\begin{equation}
\int d \mu(z) \,e^{\xi \bar{z}} z^n=\xi^n.
\end{equation}
By generalization of this identity to two complex variables we find expansion (\ref{eqn1}).

In the present paper, we are going not only proof the identity (\ref{eqn1}), but also derive the $q$-analogue of this identity.

Recently, \cite{Pashaev} we have introduced a complex function $f(z;q)$ of one complex variable $z$ according to equation
\begin{equation} D_q^{\bar{z}} f(z;q)= \frac{1}{2} \left( D_q^x+i D_{\frac{1}{q}}^y \right)f(z;q)=0 \end{equation}
and we called it as  the  \textit{$q$-analytic function}. We have described a wide class of these functions, which are not analytic in the usual sense (\ref{eqn0}).

\textbf{Example: }Complex $q$-binomial $$(x+iy)_q^2 = (x+iy)(x+iqy)= \frac{1}{2}\left((1+q)z^2+(1-q)\bar{z}z\right)$$ is not analytic $\frac{\partial}{\partial \bar{z}} (x+iy)_q^2\neq 0$, but $q$-analytic $D_q^{\bar{z}}(x+iy)_q^2=0.$\\

 A complex function $f(z,w)$ of two complex variables $z$ and $w$, analytic in both variables  and satisfying equation
 \begin{equation}\bar{D}_{z,w} f(z,w)= \frac{1}{2} \left( D_q^z+iD_{\frac{1}{q}}^w\right) f(z,w)=0\label{doubleqanaltyic}\end{equation}
we call the \textit{double $q$-analytic} function.

\textbf{Example: }Complex $q$-binomial
$$(z+iw)_q^2=z^2 +[2]_qiwz-qw^2$$ is analytic in $z$ and $w$, since $\frac{\partial}{\partial \bar{z}} (z+iw)_q^2=\frac{\partial}{\partial \bar{w}} (z+iw)_q^2 =0$ and double $q$-analytic due to $\bar{D}_{z,w}(z+iw)_q^2 =0.$

We will show that complex $q$-binomial $(z+iw)_q^n,$  for $n$-positive integer, is double $q$-analytic. This is why
 any convergent power series
 $$f(z,w)_q=\sum_{n=0}^\infty a_n (z+iw)_q^n$$ represents a double $q$-analytic function.
As a central result of the paper, we find an expansion of  the $q$-binomial in terms of  $q$-Hermite polynomials

\begin{equation}
(z+i w)_q^n = \frac{1}{[2]_q^n} \sum_{k=0}^n  {n \brack k}_q i^k q^{\frac{k(k-1)}{2}} H_{n-k}(z;q) H_k(q w,\frac{1}{q})\label{eqn2}.
\end{equation}
In the limit $q\rightarrow 1,$ this formula reduces to the Hermite binomial formula (\ref{eqn1}).

As an application of our results, we consider $q$-traveling waves in the form of power series of real $q$-binomial $(x\pm ct)_q^n.$
For these traveling waves we find $q$-Hermite binomial expansion in terms of $x$ and $t$
variables. We notice that in contrast to usual traveling waves, the $q$-traveling waves are not preserving the shape during evolution.
The q-traveling waves are subject to  $q$-wave equation
$$\left( (D_{\frac{1}{q}}^t)^2- c^2( D_q^x)^2\right) u(x,t)=0,$$
for which we solve IVP in the $q$-D'Alembert form with the Jackson integral representation.

The paper organized as follows. In Section 2 we derive the  Hermite binomial formula (\ref{eqn1}). In Section 3 we introduce $q$-Hermite polynomials   and derive the $q$-Hermite binomial formula.
The $q$-binomials are the
double $q$-analytic functions, as we demonstrate in Section 4. In Section 5 we discuss relation between  the double $q$-analytic binomial and  the   $q$-analytic one. In Section 6 we introduce $q$-traveling waves in terms of real $q$-binomials and solve IVP for $q$-wave equation in D'Alembert form. We illustrate our results for different initial values by particular examples   in Section 7. Finally, in Section 8, we derive $q$-Hermite polynomial expansion of $q$-traveling waves.

\section{Hermite Binomial Formula}
In this section, we derive the Hermite binomial formula (\ref{eqn1}). We start with the following Lemma:
\begin{lem}
For arbitrary complex numbers $ \xi$  and $\eta $:
\begin{equation} \sum_{k=0}^n { n \choose k} (-1)^k \, \frac{d^{n-k}}{d \xi^{n-k}}   \frac{d^k}{d \eta^k} e^{\xi^2/4-\eta^2/4}= \left(\frac{\xi+\eta}{2}\right)^n e^{\xi^2/4-\eta^2/4} \label{4}
\end{equation}
\end{lem}
\begin{prf}
\begin{equation} \sum_{k=0}^n { n \choose k} (-1)^k \, \frac{d^{n-k}}{d \xi^{n-k}}   \frac{d^k}{d \eta^k} e^{\xi^2/4-\eta^2/4}=\left( \frac{d}{d\xi}- \frac{d}{d\eta}\right)^n e^{\xi^2/4-\eta^2/4}
\end{equation}
By changing variables $\xi$ and $\eta$ to $\lambda$ and $\mu$ according to $$\lambda+\mu=\xi,\,\,\,\,\,\,\,\, \lambda-\mu=\eta, $$ we have
$$\lambda \, \mu=\frac{\xi^2-\eta^2}{4}= (\frac{\xi-\eta}{2})  (\frac{\xi+\eta}{2}) $$
and
\begin{equation} \left( \frac{d}{d\xi}- \frac{d}{d\eta}\right)^n e^{\lambda \mu}= \left( \frac{d}{d \mu}\right)^n e^{\lambda \mu }= \lambda^n e^{\lambda \mu }= \left(\frac{\xi+\eta}{2}\right)^n e^{\xi^2/4-\eta^2/4}.\,\,\,\,\,\,\,\,\,\,\,\,\,\,\, \diamondsuit  \end{equation}
\end{prf}

\begin{cor}
\begin{equation} \sum_{k=0}^n { n \choose k} (-1)^k \, \left(\frac{d^{n-k}} {d \xi^{n-k}}e^{\xi^2/4}\right)  \left( \frac{d^k}{d \xi^k} e^{-\xi^2/4}\right)= \xi^n \end{equation}
\end{cor}

\begin{prf}
By taking the limit of (\ref{4}) as $\eta\rightarrow \xi\Rightarrow \mu\rightarrow 0$ and $ \lambda \rightarrow \xi,$
\begin{equation}\lim_{\eta\rightarrow \xi}  \sum_{k=0}^n { n \choose k} (-1)^k \, \frac{d^{n-k}}{d \xi^{n-k}}   \frac{d^k}{d \eta^k} e^{\xi^2/4-\eta^2/4}= \lim_{\mu \rightarrow 0, \lambda \rightarrow \xi} \lambda ^n e^{\lambda \mu}=\xi^n.
\end{equation}
\end{prf}
\begin{lem}
For an arbitrary complex number $\xi$ and  $n=1,2,...,$,
\begin{equation} \int dz d\bar{z} e^{-z \bar{z}} e^{\xi \bar{z}} z^n= \xi^n. \label{2}
\end{equation}
\end{lem}

\begin{prf}
By changing complex coordinates to  the Cartesian ones,  the integral is expressed in terms of the sum
\begin{eqnarray}\int dz d\bar{z} e^{-z \bar{z}} e^{\xi \bar{z}} z^n&=& \frac{1}{\pi} \int dx dy e^{\xi (x-iy)}e^{-(x^2+y^2)}(x+iy)^n  \nonumber \\
& =& \frac{1}{\pi} \sum_{k=0}^n { n \choose k} i^k \int dx\,\,\, x^{n-k} e^{-x^2} e^{\xi x} \int dy\,\,\, y^k e^{-y^2} e^{-i\xi y} \nonumber \\
& =& \frac{1}{\pi} \sum_{k=0}^n { n \choose k} i^k \, \frac{d^{n-k}}{d \xi^{n-k}} \int dx\,\,\, e^{-x^2+\xi x} \frac{d^k}{d(-i \xi)^k} \int dy\,\,\, e^{-y^2-i\xi y} \nonumber \\
& =&  \sum_{k=0}^n { n \choose k} (-1)^k \, \frac{d^{n-k}}{d \xi^{n-k}}  e^{\xi^2/4}  \frac{d^k}{d \xi^k} e^{-\xi^2/4}\label{1},
\end{eqnarray}
where we have used the Gaussian integrals  $$\int^\infty_{-\infty} e^{-x^2+a x} dx= \sqrt{\pi} e^{a^2/4}$$ and $$\int^\infty_{-\infty} e^{-x^2+ibx} dx= \sqrt{\pi} e^{-b^2/4}.$$
By using Corollary 2.0.3. we find desired result (\ref{2}).

\end{prf}
We can generalize this result for an arbitrary analytic function given by power series $f(z)=\sum_{n=0}^\infty a_n z^n$,
so that,
 $$\int d \mu(z) e^{\xi \bar{z}} f(z)=f(\xi).$$

The above proof implies some interesting binomial identity for Hermite polynomials.
We start from the Rodrigues formula for Hermite polynomials:

\begin{equation} H_n(z)=(-1)^n e^{z^2} \frac{d^n}{d z^n} e^{-z^2}\label{hermite1}\end{equation} and by replacing $z\rightarrow iz$
\begin{equation}H_n(i z)=i^n e^{-z^2} \frac{d^n}{d z^n} e^{z^2}. \label{hermite2} \end{equation}
Then we have the following
\begin{iden}

\begin{equation}\frac{1}{2^n} \sum^n_{k=0} { n \choose k}(-i)^{n-k} H_{n-k}(\frac{i}{2} \xi) H_k(\frac{\xi}{2})=\xi^n. \label{5} \end{equation}
 \end{iden}
 \begin{prf}
According to previous proof
\begin{equation}\xi^n= \int dz d\bar{z} e^{-z \bar{z}} e^{\xi \bar{z}} z^n=\sum_{k=0}^n { n \choose k} (-1)^k \, \frac{d^{n-k}}{d \xi^{n-k}}  e^{\xi^2/4}  \frac{d^k}{d \xi^k} e^{-\xi^2/4}. \end{equation}
By inserting  $1 = e^{\xi^2/4} e^{-\xi^2/4}$ and using Rodrigues formulas we have
\begin{eqnarray} & & \sum_{k=0}^n { n \choose k} (-1)^k \,e^{\xi^2/4}  \left(e^{-\xi^2/4}\frac{d^{n-k}}{d \xi^{n-k}}  e^{\xi^2/4}\right) \left(  \frac{d^k}{d \xi^k} e^{-\xi^2/4}\right) \nonumber \\= & & \sum_{k=0}^n { n \choose k} (-1)^k \, \left(\frac{1}{2i}\right)^{n-k} H_{n-k} \left(\frac{i}{2} \xi\right)\left( e^{\xi^2/4} \frac{d^k}{d \xi^k} e^{-\xi^2/4}\right) \nonumber \\= & & \sum_{k=0}^n { n \choose k} (-1)^k \, \left(\frac{1}{2i}\right)^{n-k} H_{n-k} \left(i \frac{\xi}{2} \right)\left(-\frac{1}{2}\right)^{k} H_{k} \left(\frac{\xi}{2}\right) \nonumber \\
=& &  \frac{1}{2^n} \sum_{k=0}^n { n \choose k} (-1)^{n-k}  i^{n-k}H_{n-k} \left(i \frac{\xi}{2}\right) H_{k} \left(\frac{\xi}{2}\right) \nonumber \\ = & & \xi^n
\,\,\,\,\,\,\,\,\,\,\,\,\,\,\,\,\,\,\,\,\,\,\,\,\,\,\,\,\,\,\,\,\,\,\,\,\,\,\,\,\,\,\,\,\,\,\,\,\,\,\,\,\,\,\diamondsuit \end{eqnarray}
\end{prf}

 Particular forms of this identity are:

\begin{equation} \frac{1}{2^{2n}} \sum^n_{k=0} { n \choose k} i^k H_{n-k}( z) H_k(-iz)=z^n, \,\,\,\,\,\, (\xi \rightarrow -2iz) \label{6}
\end{equation}
and (by reductions $\xi=x$ and $\xi= iy$),
\begin{equation}\frac{1}{2^n} \sum^n_{k=0} { n \choose k} i^{n-k}H_{n-k}\left( \frac{ix}{2}\right) H_k\left(\frac{x}{2}\right)=x^n \label{7},
\end{equation}
\begin{equation}\frac{1}{2^n} \sum^n_{k=0} { n \choose k} i^{n-k}H_{n-k}\left( \frac{-y}{2}\right) H_k\left(i\frac{y}{2}\right)=i^n y^n \label{8}.
\end{equation}


This result can be generalized to the Hermite binomial formula with two complex variables $z$ and $w$

\begin{iden}
\begin{equation}(z+ i w )^n=\frac{1}{2^n} \sum^n_{k=0} { n \choose k }i^k H_{n-k}(z) H_k(w) \label{hermiteformula} \end{equation}
\end{iden}

\begin{prf}
From generating function for Hermite polynomials
\begin{equation} g(z,t)= e^{-t^2+2t z}= \sum_{n=0}^\infty H_n(z) \frac{t^n}{n!} \label{gen1}
\end{equation}
and
\begin{equation} g(w,\tau)= e^{-\tau^2+2\tau w}= \sum_{k=0}^\infty H_k(w) \frac{\tau^k}{k!} \label{gen2}\end{equation}
by changing variable $\tau=i t$ in the second one and multiplying (\ref{gen1}) and (\ref{gen2}) we have
\begin{eqnarray} g(z,t)\, g(w,it)= e^{2 t(z+i w)}=\sum_{l=0}^\infty \sum_{k=0}^\infty \frac{H_l(z)}{l!} \frac{H_k(w)i^k}{k!} t^{l+k}.
\end{eqnarray}
By  changing the order of double sum with
 $l+k=n,$ and expanding the left hand side in $t$ we get
\begin{equation} g(z,t)\, g(w,it)= \sum_{n=0}^\infty \frac{2^n t^n (z+i w)}{n!} =\sum_{n=0}^\infty \frac{t^n}{n!} \sum_{k=0}^n { n \choose k} H_{n-k}(z) H_k(w) i^k.
\end{equation}
By equating terms of the same power $t^n$ we obtain the desired result (\ref{hermiteformula}):
$$(z+ i w )^n=\frac{1}{2^n} \sum^n_{k=0} { n \choose k }i^k H_{n-k}(z) H_k(w).\,\,\,\,\,\,\,\,\,\,\, \diamondsuit $$
\end{prf}

Here we can give also another proof of this relation by using the holomorphic Laplace equation.
\begin{prf}
$\zeta^n \equiv(z+ i w)^n$ is double analytic function of two complex variables $z$ and $w$. Therefore, due to (\ref{2analytic}) it satisfies the holomorphic  Laplace equation $\Delta \zeta^n=0,$
where $\Delta \equiv \frac{\partial^2}{\partial z^2}+ \frac{\partial^2}{\partial w^2} = \left(\frac{\partial}{\partial z}- \frac{\partial}{\partial w}\right)
\left(\frac{\partial}{\partial z}+ \frac{\partial}{\partial w}\right)
$,
which implies $\Delta ^k \zeta^n=0$
for  every $k=1,2...$
Expanding
\begin{equation}
 e^{-\frac{1}{4}\Delta}(z+i w)^n= \left( 1-\frac{1}{4} \Delta + \frac{(-\frac{1}{4}\Delta)^2}{2!}+...+ \frac{(-\frac{1}{4}\Delta)^n}{n!}+...\right)(z+i w)^n = (z+i w)^n \end{equation}
we have
\begin{eqnarray}
e^{-\frac{1}{4}\Delta}(z+i w)^n&= &e^{-\frac{1}{4}\frac{d^2}{dz^2}}\,e^{-\frac{1}{4}\frac{d^2}{dw^2}} \sum_{k=0}^n { n \choose k } z^{n-k} i^k w^k
\nonumber \\
&=& \sum_{k=0}^n { n \choose k } i^k\, \left( e^{-\frac{1}{4}\frac{d^2}{dz^2}}z^{n-k}\right) \,\left( e^{-\frac{1}{4}\frac{d^2}{dw^2}}w^k \right).
\end{eqnarray}

Then due to relation:
$$H_n(z)= 2^n e^{-\frac{1}{4}\frac{d^2}{dz^2}} z^n$$
we get the Hermit binomial formula
\begin{equation}
(z+i w)^n= \frac{1}{2^n} \sum_{k=0}^n { n \choose k } H_{n-k}(z) H_k(w) i^k.
\end{equation}
\end{prf}
\section{q-Hermite Binomial Formula}
Here we are going to generalize the above formula to the q-binomial case. For this, first we need to introduce the q-Hermite polynomials.
\subsection{q-Hermite Polynomials}
In paper \cite{Nalci}, studying the q-heat and q-Burgers equations, we have defined the $q$-Hermite polynomials according to generating function

\begin{equation} e_q(-t^2) e_q([2]_q t x) = \sum^{\infty}_{n=0} H_n(x;q) \frac{t^n}{[n]_q!},\label{genfunc}
 \end{equation}
 where $$e_q(x)= \sum_{n=0}^\infty \frac{x^n}{[n]_q!},\,\,\,\,\,\,\,\,E_q(x)= \sum_{n=0}^\infty q^{\frac{n(n-1)}{2}}\frac{x^n}{[n]_q!} $$
are Jackson's $q$-exponential functions and $q$-numbers and $q$-factorials are defined as follows: $$[n]_q=\frac{q^n-1}{q-1},\,\,\,\,\,\, [n]_q!= [1]_q [2]_q...[n]_q.$$

By q-differentiating the generating function (\ref{genfunc})  according to $x$ and $t$ we have the recurrence relations,
correspondingly

\begin{equation} D_x H_n (x;q) = [2]_q [n]_q H_{n-1}(x;q),\label{rec1}
\end{equation}
\begin{equation}
H_{n+1} (x;q) = [2]_q \,x  H_{n}(x;q) - [n]_q \,H_{n-1}(q x;q)\nonumber \\
- [n]_q \,q^{\frac{n+1}{2}} H_{n-1}(\sqrt{q} x;q).\label{rec2}
\end{equation}

We get also  the special values
\begin{eqnarray} &&H_{2n}(0;q) = (-1)^n \frac{[2n]_q!}{[n]_q!},\\
&&H_{2n+1}(0;q) = 0, \end{eqnarray}
and the parity relations
\begin{equation} H_n (-x;q) = (-1)^n H_n (x;q). \end{equation}
 For more details we refer to paper \cite{Nalci}.
 First few polynomials are
\begin{eqnarray} H_0 (x;q) = 1, \,\,\,\,\,H_1 (x;q) = [2]_q \,x \nonumber,\end{eqnarray}
\begin{eqnarray} H_2(x;q) = [2]^2_q \, x^2 - [2]_q,\,\,\,\,\,H_3(x;q) = [2]^3_q \,x^3 - [2]^2_q [3]_q \,x\nonumber, \end{eqnarray}
\begin{eqnarray} H_4(x;q) = [2]^4_q \,x^4 - [2]^2_q [3]_q [4]_q\,x^2 + [2]_q [3]_q [2]_{q^2} \nonumber ,\end{eqnarray}
and when $q\rightarrow 1$ these polynomials reduce to the standard Hermite polynomials.

The generating function (\ref{genfunc}) for $t=1$ gives expansion of $q$-exponential function in terms of $q$-Hermite polynomials

\begin{equation} \sum_{n=0}^\infty \frac{H_n(x;q)}{[n]_q!} = \frac{e_q([2]_q x)}{E_q(1)}.
\end{equation}
In the limiting case $q \rightarrow 1$ it gives expansion of exponential
$$ e^{x}= e \sum_{n=0}^\infty \frac{H_n\left(\frac{x}{2}\right)}{n!}$$
and for $x=1$ we find  for  Euler's number $e$ as:
$$e=\sum_{n=0}^\infty \frac{H_n(1)}{n!}.$$
 For the Jackson q-exponential function we have the q-analog of this expansion:
$$  e_q(x)  =   e_{\frac{1}{q}} \sum_{n=0}^\infty \frac{H_n\left(\frac{x}{[2]_q};q\right)}{[n]_q!},       $$
and as follows
\begin{equation}  \frac{e_q}{e_{\frac{1}{q}}} = e_q(1) e_{q}(-1) =   \sum_{n=0}^\infty \frac{H_n\left(\frac{1}{[2]_q};q\right)}{[n]_q!},  \label{qEuler}    \end{equation}
where
$$ e_q = \sum_{n=0}^\infty \frac{1}{[n]_q!}$$
is the q-analog of Euler's number $e$.
Relation (\ref{qEuler}) should be compared with the next one
$$ e_q(1) e_{q}(-1) = e_{q^2} \left( \frac{1-q}{1+q}\right) = \prod^\infty_{k=1} \frac{1}{1 - (1-q^2)q^k}$$
which is coming for  $x=1$ from the following identity.

\begin{iden}

$$ e_q(x)e_q(-x) = e_{q^2} \left( \frac{1-q}{1+q}x^2\right).$$
\end{iden}
\begin{prf}
We expand and change order of the double sum
 \begin{eqnarray}e_q(x)e_q(-x)  = \sum_{k=0}^\infty \frac{x^k}{[k]_q!}\sum_{l=0}^\infty \frac{(-1)^l x^l}{[l]_q!} = \sum_{k=0}^\infty \sum_{l=0}^\infty (-1)^l \frac{x^{k+l}}{[k]_q! [l]_q!} =\\
 \sum_{n=0}^\infty   \frac{x^n}{[n]_q!}  \sum_{k=0}^n  (-1)^k \frac{[n]_q!}{[k]_q! [n-k]_q!} = \sum_{n=0}^\infty   \frac{x^n}{[n]_q!}  \sum_{k=0}^n  {n \brack k}_q (-1)^k .\end{eqnarray}
 Splitting the first sum to the even $n=2m$ and to the odd parts $n=2m+1$ we have
 \begin{eqnarray}   \sum_{m=0}^\infty   \frac{x^{2m}}{[2m]_q!}  \sum_{k=0}^{2m}  {2m \brack k}_q (-1)^k  +
 \sum_{m=0}^\infty   \frac{x^{2m+1}}{[2m+1]_q!}  \sum_{k=0}^{2m+1}  {2m+1 \brack k}_q (-1)^k.\end{eqnarray}
 Due to known identities \cite{Kac et al.}
 $$ \sum_{k=0}^{2m+1}  {2m+1 \brack k}_q (-1)^k=0,     $$
 $$\sum_{k=0}^{2m}  {2m \brack k}_q (-1)^k= (1-q) (1-q^3)...(1-q^{2m-3}) (1-q^{2m-1}),$$
the second sum vanishes and for the first sum we get
  $$ \sum_{m=0}^\infty \frac{x^{2m}}{[2m]_q!} (1-q) (1-q^3)...(1-q^{2m-3}).  $$
 Since $[2 m]_q= [2]_q [m]_{q^2}$, we can rewrite  this sum as

  $$\sum_{m=0}^\infty \frac{x^{2m} (1-q)^m}{[2]_q^m \,\,[m]_{q^2}! } = e_{q^2} \left( \frac{(1-q)x^2}{[2]_q}\right).$$
  Finally, $$ e_q(x)e_q(-x) = e_{q^2} \left( \frac{1-q}{1+q}x^2\right).$$
\end{prf}

\subsection{$q$-Hermite Binomials}
Here we formulate our main result as the $q$-Hermite binomial identity.
\begin{iden}
The $q$-analogue of identity (\ref{hermiteformula}), giving $q$-binomial expansion in terms of $q$-Hermite polynomials is
\begin{equation}
(z+iw)_q^n = \frac{1}{[2]_q^n} \sum_{k=0}^n  {n \brack k}_q i^k q^{\frac{k(k-1)}{2}} H_{n-k}(z;q) H_k(qw,\frac{1}{q}) \label{qhermiteformula}
\end{equation}
\end{iden}

\begin{prf}
By using the generating function for q-Hermite polynomials (\ref{genfunc}) and replacing $x\rightarrow Z$ we obtain
\begin{equation} e_q(-t^2) e_q([2]_q Z t)=\sum_{n=0}^\infty H_n(Z;q) \frac{t^n}{[n]_q!}. \label{genhermite} \end{equation}
In this formula we replace
 $t\rightarrow i t ,$ $ Z\rightarrow W,$ $q\rightarrow 1/q$ so that
\begin{equation} e_{\frac{1}{q}}(t^2) e_{\frac{1}{q}}([2]_{\frac{1}{q}} i W t)=\sum_{n=0}^\infty H_n(W;{\frac{1}{q}}) i^n \frac{t^n}{[n]_{\frac{1}{q}}!}. \label{genhermite2} \end{equation}
Multiplying (\ref{genhermite}) with (\ref{genhermite2}) and using  factorization of $q$-exponential functions

\begin{equation}e_q(x) e_{\frac{1}{q}}(y)=  \sum_{n=0}^\infty \frac{(x+y)_q^n}{[n]_q!}\equiv e_q(x+y)_q \label{expfactor}\end{equation}
leading to
$$e_q(-t^2) e_{\frac{1}{q}}(t^2)=e_q(0)_q=1, $$
we find

\begin{equation}
e_{q}( t ([2]_q Z+ [2]_{\frac{1}{q}}i W) )_q = \sum_{l=0}^\infty \sum_{k=0}^\infty \frac{H_l(Z;q)}{[l]_q!} \frac{H_k(W;\frac{1}{q} )i^k}{[k]_{\frac{1}{q}}!} t^{l+k}.
 \end{equation}
 By changing order of the double sum and expanding the left hand side in $t,$ we get
\begin{equation}
\sum_{n=0}^\infty \frac{t^n
( [2]_q Z+ [2]_{\frac{1}{q}}i W)_q^n}{[n]_q!}= \sum_{n=0}^\infty \frac{t^n}{[n]_q!} \sum_{k=0}^n {n \brack k}_q q^{\frac{k(k-1)}{2}}  H_{n-k}(Z;q) H_k(W;\frac{1}{q}) i^k.
\end{equation}
Then, at power $t^n$ we have identity

\begin{equation}
\left(Z+i\frac{W}{q }\right)_q^n= \frac{1}{[2]_q^n} \sum_{k=0}^n{n \brack k}_q q^{\frac{k(k-1)}{2}}  H_{n-k}(Z;q) H_k(W;\frac{1}{q}) i^k,
\end{equation}
where $$[k]_{\frac{1}{q}}= \frac{1}{q^{k-1}} [k]_q,\,\,\,\,\,\,\,\,[k]_{\frac{1}{q}}!= \frac{1}{q^{\frac{k(k-1)}{2}}} [k]_q!$$
By replacing $Z=z$ and $\frac{W}{q}=w$ the desired result is obtained
$$(z+iw)_q^n = \frac{1}{[2]_q^n} \sum_{k=0}^n{n \brack k}_q q^{\frac{k(k-1)}{2}}  H_{n-k}(z;q) H_k(qw;\frac{1}{q}) i^k.\,\,\,\,\,\,\,\,\,\,\, \diamondsuit $$
\end{prf}

\section{Double q-Analytic Function}

Here we consider a class of complex valued functions of two complex variables, $z$ and $w$, (or four real variables), analytic in these variables
$\frac{\partial}{\partial\bar z} f = \frac{\partial}{\partial\bar w} f = 0$.

\begin{defn}
A complex-valued function $f(z,w)$ of four real variables is called the double analytic in a region if the following identity holds in the region:
\begin{equation}
\bar{\partial}_{z,w}f \equiv \frac{1}{2}(\partial_{z} +i \partial_w)f=0 ,
\end{equation}
where $$\partial_{z} f= \frac{1}{2}(\partial_{x} - i \partial_{y}),\,\,\,\,\partial_{w} f= \frac{1}{2}(\partial_{u} - i \partial_{v})$$
and $z=x+iy,\,\,\,\,\,w=u+iv.$
\end{defn}

\begin{defn}
A complex-valued function $f(z,w)$ of four real variables is called the double $q$-analytic in a region if the following identity holds in the region:
\begin{equation}
\bar{D}_{z,w}f= \frac{1}{2}(D_q^{z} +i D_{\frac{1}{q}}^w)f=0 ,
\end{equation}
where $$D_q^{z} f(z, w)= \frac{f(qz, w)-f(z, w)}{(q-1)z},\,\,\,\,\,\,\,\,D_{\frac{1}{q}}^{w} f(z, w)= \frac{f(z, \frac{w}{q})-f(z, w)}{(\frac{1}{q}-1)w}$$
and $z=x+iy,\,\,\,\,\,w=u+iv.$
\end{defn}
Here we should notice that
$$ D_q^{z} \neq \frac{1}{2}(D_q^{x}  - i   D_q^{y}),\,\,\,\,D_q^{w} \neq \frac{1}{2}(D_q^{u}  - i   D_q^{v}).$$

The simplest set of double $q$-analytic functions is given by complex $q$-binomials
$$(z+ i w)_q^n \equiv (z+ i w)(z+i q w)(z+iq^2 w)...(z+i q^{n-1} w) = \sum^n_{k=0} {n \brack k }_q q^{k(k-1)/2} i^k z^{n-k} w^k$$
satisfying
 $$\frac{1}{2}(D_q^{z} +i D_{\frac{1}{q}}^w)(z+ i w)_q^n=0 ,$$
 and
$$\frac{1}{2}(D_q^{z} -i D_{\frac{1}{q}}^w)(z+ i w)_q^n=[n]_q(z+iw)_q^{n-1}.$$
\begin{prf}

\begin{eqnarray} && \frac{1}{2}(D_q^{z} -i D_{\frac{1}{q}}^w)(z+ i w)_q^n \nonumber \\ & = & \frac{1}{2}\left( \sum^n_{k=0} {n \brack k }_q q^{k(k-1)/2} i^k \left(D_q^z\,\, z^{n-k} \right)w^k-i
\sum^n_{k=0} {n \brack k }_q q^{k(k-1)/2} i^k z^{n-k} \left(D_{\frac{1}{q}}^w\,\, w^k  \right)\right) \nonumber \\
&=& \frac{1}{2} \left( \sum^{n-1}_{k=0} {n \brack k }_q  q^{k(k-1)/2} [n-k]_q z^{n-k-1} i^k w^k -i \sum^n_{k'=1} {n \brack k' }_q q^{{k'}({k'}-1)/2} i^{k'} z^{n-{k'}} \frac{w^{{k'}-1}}{q^{{k'}-1}}  [k']_q \right) \nonumber \\
&=&\frac{1}{2}\left( \sum^{n-1}_{k=0} {n \brack k }_q  q^{k(k-1)/2} [n-k]_q z^{n-k-1} i^k w^k -i \sum^{n-1}_{k=0} {n \brack k+1 }_q q^{k(k+1)/2} i^{k+1} z^{n-k-1} \frac{w^{k}}{q^k}  [k+1]_q \right)\nonumber \\
&=& \frac{1}{2}\sum^{n-1}_{k=0} \left( {n \brack k }_q q^{k(k-1)/2}[n-k]_q -i {n \brack k+1 }_q q^{k(k+1)/2} i \frac{1}{q^k}[k+1]_q\right) z^{n-k-1} i^k w^k \nonumber \\
&=& \frac{1}{2}  \sum^{n-1}_{k=0} 2 \left( \frac{[n]!}{[n-k-1]! [k]!} q^{k(k-1)/2} \right)z^{n-k-1} i^k w^k  \nonumber \\
&=& [n] \sum^{n-1}_{k=0}{n-1 \brack k }_q q^{k(k-1)/2}z^{n-k-1} i^k w^k  \nonumber \\
&=& [n](z+i w)_q^{n-1}\nonumber.
\end{eqnarray}
\end{prf}

From above result follows that
any convergent power series $$f(z+iw)_q=\sum_{n=0}^\infty a_n (z+iw)_q^n$$
determines a double $q$-analytic function.
Since our relation (\ref{qhermiteformula}) shows expansion of double $q$-analytic  $q$-binomials in terms of $q$-Hermite polynomials,
 it also gives expansion of any double $q$-analytic function in terms of the analytic polynomials.

\textbf{Examples:}
For $n=1:$
\begin{eqnarray}
(z+iw)_q^1 = z+ i w= \frac{1}{[2]_q} \sum_{k=0}^1 {1 \brack k}_q q^{\frac{k(k-1)}{2}} H_{1-k}(z;q) H_k(qw;\frac{1}{q}) i^k  \nonumber  \\
=  \frac{1}{[2]_q} \left({1 \brack 0}_q  H_1(z;q) H_0(qw;\frac{1}{q}) + {1 \brack 1}_q   H_0(z;q) H_1(qw;\frac{1}{q})i  \right) \nonumber \\
.
\end{eqnarray}

For $n=2:$
\begin{eqnarray}(z+iw)_q^2 & =& (z+iw) (z+i qw)= z^2+i[2]_q z w -q w^2\nonumber \\ & =& \frac{1}{[2]_q^2} \sum_{k=0}^2 {2 \brack k}_q q^{\frac{k(k-1)}{2}} H_{2-k}(z;q) H_k(qw;\frac{1}{q}) i^k \nonumber \\
 &=&\frac{1}{[2]_q^2} \left({2 \brack 0}_q  H_2(z;q) H_0(qw;\frac{1}{q}) + {2 \brack 1}_q   H_1(z;q) H_1(qw;\frac{1}{q})i+ {2 \brack 2}_q  q H_0(z;q) H_2(qw;\frac{1}{q}) i^2\right) \nonumber \end{eqnarray}

\subsection{$q$-Holomorphic Laplacian}

Another proof of identity (\ref{qhermiteformula}) can be done by noticing that  $q$-binomial $(z+i w)_q^n$ is double $q$-analytic function.
Then we can use the following identity and complex q-Laplace equation.
\begin{iden}
\begin{equation}
e_q (-\frac{1}{[2]_q^2} \Delta_q )_q (z+iw)_q^n= (z+iw)_q^n \label{identity}
\end{equation}
\end{iden}
\begin{prf}

Since $(z+iw)_q^n$ is $q$-analytic function, it satisfies the $q$-Laplace equation $$\Delta_q (z+iw)_q^n=0$$
and
\begin{eqnarray}
e_q (-\frac{1}{[2]_q^2} \Delta_q )_q (z+iw)_q^n &=& e_q (-\frac{1}{[2]_q^2} \left( (D_q^z)^2+ (D_{\frac{1}{q}}^w)^2\right)) _q (z+iw)_q^n\nonumber \\
&=& \sum_{n=0}^\infty \frac{1}{[n]_q!} \left(-\frac{1}{[2]_q^2}\right)^n \left(  (D_q^z)^2 + (D_{\frac{1}{q}}^w)^2\right)_q^n (z+iw)_q^n\nonumber \\
&=&  \sum_{n=0}^\infty \frac{1}{[n]_q!}\left(-\frac{1}{[2]_q^2}\right)^n \left( \Delta_q \right)_q^n(z+iw)_q^n,
\end{eqnarray}
where $\left( \Delta_q \right)_q^n = \Delta_q \cdot \Delta_q^{(1)} \cdot \Delta_q^{(2)} \cdot ... \cdot \Delta_q^{(n-1)} $
and $$\Delta_q=(D_q^z)^2+ (D_{\frac{1}{q}}^w)^2,\,\,\,\,\,\,\,\,\,\,\Delta_q^{(1)}= (D_q^z)^2+q (D_{\frac{1}{q}}^w)^2,\,\,\,...,\,\, \Delta_q^{(n-1)}= (D_q^z)^2+q^{n-1} (D_{\frac{1}{q}}^w)^2.   $$
Using the fact that $\left( \Delta_q\right)_q^m (z+iw)_q^n=0, \forall m=1,2,...$,
only the first term in expansion survives,
then we get desired result.
\end{prf}
Due to (\ref{expfactor}) we can factorize q-exponential operator function as
\begin{eqnarray}
& &e_q \left(-\frac{1}{[2]_q^2} \Delta_q \right)_q (z+iw)_q^n \nonumber \\
&=& e_q \left( -\frac{1}{[2]_q^2} (D_q^z)^2 \right) e_{\frac{1}{q}} \left( -\frac{1}{[2]_q^2} (D_{\frac{1}{q}}^w)^2 \right) (z+iw)_q^n \nonumber \\
&=& \sum_{k=0}^n {n \brack k}_q q^{\frac{k(k-1)}{2}} i^k e_q \left( -\frac{1}{[2]_q^2} (D_q^z)^2 \right)z^{n-k} e_{\frac{1}{q}} \left( -\frac{1}{[2]_q^2} (D_{\frac{1}{q}}^w)^2 \right)w^k. \label{a}
\end{eqnarray}

By using the generating function of q-Hermite Polynomials (\ref{genhermite}) we have the following identity:
\begin{equation}
H_n(x;q)= [2]_q^n e_q \left( -\frac{1}{[2]_q^2} (D_q^x)^2\right) x^n,
\end{equation}
which gives
\begin{equation}
e_q \left( -\frac{1}{[2]_q^2} (D_q^z)^2 \right)z^{n-k}= \frac{1}{[2]_q^{n-k}} H_{n-k} (z;q),
\end{equation}
and
\begin{equation}
e_{\frac{1}{q}} \left( -\frac{1}{[2]_q^2} (D_{\frac{1}{q}}^w)^2 \right)w^k= \frac{1}{[2]_{\frac{1}{q}}^{k}q^k} H_{k} (qw;\frac{1}{q}),
\end{equation}
where $D_{\frac{1}{q}}^{qw}= \frac{1}{q} D_{\frac{1}{q}}^w.$
Substituting into  (\ref{a}), we get
$$e_q (-\frac{1}{[2]_q^2} \Delta_q )_q (z+iw)_q^n= \sum_{k=0}^n {n \brack k}_q q^{\frac{k(k-1)}{2}} i^k \frac{1}{[2]_q^{n-k}} H_{n-k} (z;q) \frac{1}{[2]_{\frac{1}{q}}^{k}q^k} H_{k} (qw;\frac{1}{q}).$$
Then, according to identity (\ref{identity})
we obtain desired result
\begin{equation}(z+iw)_q^n= \frac{1}{[2]_q^n}\sum_{k=0}^n {n \brack k}_q q^{\frac{k(k-1)}{2}}H_{n-k} (z;q) H_{k} (qw;\frac{1}{q}) i^k.\,\,\,\,\,\,\,\,\,\,\,\diamondsuit  \label{formula} \end{equation}

As a particular case of our binomial formula, we can find $q$-Hermite binomial expansion for the $q$-analytic binomial $(x+iy)^n$ as well. If in (\ref{qhermiteformula}) we   replace $z\rightarrow x$ and $w \rightarrow y,$ then we get
\begin{equation}(x+iy)_q^n= \frac{1}{[2]_q^n}\sum_{k=0}^n {n \brack k}_q q^{\frac{k(k-1)}{2}}H_{n-k} (x;q) H_{k} (qy;\frac{1}{q}) i^k.  \end{equation}
Since a $q$-analytic function is determined by power series in $q$-binomials \cite{Pashaev},
this formula allows us to get expansion of an arbitrary $q$-analytic function in terms of real $q$-Hermite polynomials.


\section{q-Traveling Waves}
As an application of $q$-binomials here we consider the $q$-analogue of traveling waves as a solution of $q$-wave equation.
\subsection{Traveling Waves:}
Real functions of two real variables
$F(x,t)=F(x \pm c t)$ called the traveling waves,  satisfy the following first order equations

$$\left(\frac{\partial}{\partial t} \mp c \frac{\partial}{\partial x}  \right)F(x \pm c t)=0. $$
It describes waves with fixed shape, prorogating with constant speed $c$ in the left and in the right direction correspondingly. The general solution of the wave equation
$$\frac{\partial^2 u}{\partial t^2}=c^2\frac{\partial^2 u}{\partial x^2}, $$
then can be written as an arbitrary superposition of these  traveling waves
$$u(x,t)= F(x+c t)+ G(x- ct).$$

\subsection{q-Traveling Waves}
Direct extension of traveling waves to q-traveling waves is not possible. This happens due to the absence in q-calculus of the chain rule and as follows, impossibility to
use moving frame as an argument of the wave function. Moreover, if we try in the Fourier harmonics $f(x,t)= e^{i (k x-\omega t)},$ replace exponential function by Jackson's $q$-exponential function $f(x,t)=e_q(i(k x-\omega t)),$ then we find that it doesn't work due to the absence of factorization for $q$-exponential function $e_q(i (k x-\omega t)) \neq  e_q(i k x)\, e_q(i \omega t).$

This is why, here we propose another way. First we observe that q-binomials
\begin{equation}(x\pm ct)_q^n= (x \pm ct)(x \pm q ct)...(x \pm q^{n-1} ct)\label{qtravbinom}\end{equation}
for $n = 0, \pm1, \pm2,... $, satisfy the first order one-directional $q$-wave equations
\begin{equation}
\left( D^t_{\frac{1}{q}} \mp c D_q^x\right) (x\pm ct)_q^n=0.\label{qtrav}
\end{equation}
Then,
the Laurent series expansion in terms of these $q$-binomials determines the q-analog of traveling waves
$$f(x\pm ct)_q = \sum_{n=-\infty}^\infty a_n (x\pm ct)_q^n.$$
Due to (\ref{qtrav}) the q-binomials (\ref{qtravbinom}) satisfy
the $q$-wave equation
 \begin{equation}\left( (D_{\frac{1}{q}}^t)^2- c^2( D_q^x)^2\right) u(x,t)=0 \label{wave}
\end{equation}
and the general solution of this equation is expressed in the form of $q$-traveling waves
\begin{equation}
u(x,t)=F(x+ct)_q+ G(x-ct)_q
\end{equation}
where
$$F(x+c t)_q= \sum_{n=-\infty}^{\infty} a_n (x+c t)_q^n$$ and $$G(x-c t)_q= \sum_{n=-\infty}^{\infty} b_n (x-c t)_q^n.$$

This allows us to solve IVP for the $q$-wave equation
\bea & & \left[\left(  D_{\frac{1}{q}}^t\right)^2 - c^2 (D_{q}^x)^2\right]u(x,t)=0, \label{ivp1} \\
& &  u(x,0)= f(x), \label{ivp2}\\
& & D_{\frac{1}{q}}^t u(x,0)=g(x), \label{ivp3}\eea where  $-\infty < x<\infty,$
in the D'Alembert form:

\bea  u(x,t)= \frac{f(x+c t)_q+f(x-c t)_q}{2}+\frac{1}{2 c} \int_{(x-ct)_q}^{(x+ct)_q} g(x') d_q x',\label{dalambert}\eea
where the Jackson integral is
\bea \int_{(x-c t)_q}^{(x+c t)_q} g(x') d_q x'&=& (1-q)(x+c t)\sum_{j=0}^\infty q^j g(q^j(x+c t))_q \nonumber \\
&-&(1-q)(x-c t)\sum_{j=0}^\infty q^j g(q^j(x-c t))_q. \eea
If the initial velocity is zero, $g(x)=0,$ the formula reduces to \bea u(x,t)= \frac{1}{2}\left(f(x+c t)_q+f(x-c t)_q\right).\eea
It should be noted here that $q$-traveling wave is not traveling wave in the standard sense and it is not preserving shape during evolution. It can be seen from simple observation. The traveling wave polynomial $(x-ct)_q^n= (x- ct)(x-q ct) (x-q^2 ct)...(x-q^{n-1}ct)$ includes the set of moving frames (as zeros of this polynomial) with re-scaled set of speeds $(c,qc,q^2c,...,q^{n-1}c).$ It means that zeros of this polynomial are moving with different speeds and therefore the shape of polynomial wave is not preserving. Only in the linear case and in the case $q=1,$ when speeds of all frames coincide, we are getting standard traveling wave.

\subsection{EXAMPLES}
In this section we are going to illustrate our results by several explicit solutions.

\textbf{Example 1:}
We consider I.V.P. for the $q$-wave equation (\ref{ivp1}) with initial functions
\bea & & u(x,0)= x^2, \nonumber \\
& & D_{\frac{1}{q}}^t u(x,0)=0 .\eea
Then the solution of the given I.V.P. for $q$-wave equation in D'Alembert form is found as
\bea  u(x,t)= x^2 +q c^2 t^2.\eea
When $q=1,$ it reduces to well-known one as superposition of two traveling wave parabolas $(x\pm ct)^2$ moving to the right and to the left with speed $c.$
Geometrically, the meaning of $q$ is the acceleration of our parabolas in vertical direction.\\

\textbf{Example 2:}
The $q$-traveling wave
\bea u(x,t)= (x-c t)_q^2&=& (x-ct) (x- q ct) \nonumber \\
&=& \left( x- \frac{[2]}{2} ct\right)^2- \frac{(q-1)^2}{4} c^2 t^2 \eea
gives solution of I.V.P. for the $q$-wave equation (\ref{ivp1}) with initial functions
\bea & & u(x,0)= x^2 \nonumber \\
& & D_{\frac{1}{q}}^t u(x,0)=-[2]_q c x .\eea

If $q=1$ in this solution we have two degenerate zeros moving with the same speed c.
In the case  $q \neq 1$,  two zeros are moving with different speeds $c$ and $q c.$ It means that, the distance between zeros is growing linearly with time as $(q-1)c t.$ The solution is the parabola,  moving in vertical direction with acceleration  $ \frac{(q-1)^2}{4} c^2,$ and in horizontal direction with constant speed  $\frac{[2]_q}{2} c.$ The area under the curve between moving zeros $x=ct$ and $x=q ct$
$$\int_{ct}^{q ct} (x-ct)_q^2 dx= - \frac{(q-1)^3 c^3}{6} t^3$$
is changing according to time as  $t^3.$

For more general initial function  $f(x)=x^n, \,\,\, n= 2,3,...$ we get $q$-traveling wave
$$u(x,t)= (x-ct)_q^n= (x-ct)(x- q ct)...(x-q^{n-1} ct)$$
 with n-zeros moving with speeds $c, q c, ..., q^{n-1}c.$ The distance between two zeros is growing as $(q^m-q^n)ct,$ and the shape of wave is changing. In parabolic case with $n=2,$ the shape of curve is not changing, but moving in horizontal direction with constant speed, and in vertical direction with constant acceleration. In contrast to this, for $n>2,$ the motion of zeros with different speeds changes the shape of the wave, and it can not be reduced to simple translation and acceleration.

\textbf{Example 3:}
Given I.V.P. for the $q$-wave equation (\ref{ivp1}) with initial functions as $q$-trigonometric functions \cite{Kac et al.}
\bea & & u(x,0)= cos_q x, \nonumber \\
& & D_{\frac{1}{q}} ^t u(x,0)=sin_q x .\eea
By using the D'Alembert form (\ref{dalambert}), after $q$-integration , we get
\bea  u(x,t)&=& \frac{1}{2} \left[cos_q(x+c t)_q+cos_q(x-c t)_q\right]+\frac{1}{2 c} \int_{(x-ct)_q}^{(x+ct)_q} sin_q(x') d_q x' \nonumber \\
&=& \frac{1}{2} \left[ \left( 1+ \frac{1}{c}\right) cos_q (x-ct)_q + \left( 1- \frac{1}{c}\right) cos_q (x+ct)_q\right]. \eea

\textbf{Example 4: $q$-Gaussian Traveling Wave}
For initial function in Gaussian form: $u(x,0)=e^{-x^2}$ in the standard case $q=1$ we have the Gaussian traveling wave $u(x,t)= e^{-(x-ct)^2}.$
For the $q$-traveling wave with Gaussian initial condition $u(x,0)=e^{-x^2},$ we have the $q$-traveling wave $$u(x,t)= \left( e^{-(x-ct)^2}\right)_q \equiv \sum_{n=0}^\infty \frac{(-1)^n}{n!} (x-ct)_q^{2 n}.$$


\subsection{q-Traveling Waves in terms of q-Hermite Polynomials}

Identity(\ref{hermiteformula})
 allows us to rewrite the traveling wave binomial in terms of Hermite polynomials as
$$(x+ct)^n=\frac{1}{2^n} \sum^n_{k=0} { n \choose k }i^k H_{n-k}(x) H_k(-ict) .$$
Its $q$-analogue for $q$-traveling wave binomial follows from (\ref{qhermiteformula})
$$(x+ct)_q^n= \frac{1}{[2]_q^n} \sum_{k=0}^n  {n \brack k}_q i^k q^{\frac{k(k-1)}{2}} H_{n-k}(x;q) H_k(-iqct,\frac{1}{q}). $$

Then, the general solution of $q$-wave equation (\ref{wave}) can be expressed in the form of $q$-Hermite polynomials

\begin{equation}
u(x,t)=F(x+ct)_q+ G(x-ct)_q,
\end{equation}
where
\begin{eqnarray}& & F(x+c t)_q= \sum_{n=-\infty}^{\infty} a_n (x+c t)_q^n= \sum_{n=-\infty}^{\infty} a_n \frac{1}{[2]_q^n} \sum_{k=0}^n  {n \brack k}_q i^k q^{\frac{k(k-1)}{2}} H_{n-k}(x;q) H_k(-iqct,\frac{1}{q}),  \nonumber \\
& & G(x-c t)_q= \sum_{n=-\infty}^{\infty} a_n (x-c t)_q^n= \sum_{n=-\infty}^{\infty} a_n  \frac{1}{[2]_q^n} \sum_{k=0}^n  {n \brack k}_q i^k q^{\frac{k(k-1)}{2}} H_{n-k}(x;q) H_k(iqct,\frac{1}{q}).\nonumber
\end{eqnarray}

It is instructive to prove the $q$-traveling wave solution
$$\left( D^t_{\frac{1}{q}} - c D_q^x\right)(x+ct)_q^n=0$$  by using $q$- Hermite binomial.
We have
\begin{eqnarray}
&&(D^t_{\frac{1}{q}} - c D_q^x)(x+ct)_q^n=\left( D^t_{\frac{1}{q}} - c D_q^x\right) \frac{1}{[2]_q^n} \sum_{k=0}^n  {n \brack k}_q i^k q^{\frac{k(k-1)}{2}} H_{n-k}(x;q) H_k(-iqct,\frac{1}{q})\nonumber \\
&&= \frac{1}{[2]_q^n} \sum_{k=0}^n  {n \brack k}_q i^k q^{\frac{k(k-1)}{2}} H_{n-k}(x;q)D_{\frac{1}{q}}^t H_k(-iqct,\frac{1}{q}) \nonumber \\
& &-c \frac{1}{[2]_q^n} \sum_{k=0}^n  {n \brack k}_q i^k q^{\frac{k(k-1)}{2}}D_q^x H_{n-k}(x;q) H_k(-iqct,\frac{1}{q}) \nonumber
\end{eqnarray}
By recursion formula for  $q$-Hermite polynomials
$$D_q^x H_n (x;q)= [2]_q [n]_q H_{n-1}(x;q)$$
we get
\begin{eqnarray}
&& (D^t_{\frac{1}{q}} - c D_q^x)(x+ct)_q^n=\frac{1}{[2]_q^n} \sum_{k'=1}^n  {n \brack k'}_q i^{k'} q^{\frac{k'(k'-1)}{2}} H_{n-k'}(x;q)[2]_{\frac{1}{q}}  [k']_{\frac{1}{q}} H_{k'-1}(-iqct,\frac{1}{q})(-iqc) \nonumber \\
&& -c \frac{1}{[2]_q^n} \sum_{k=0}^{n-1}  {n \brack k}_q i^k q^{\frac{k(k-1)}{2}} [2]_q [n-k]_q H_{n-k-1}(x;q) H_k(-iqct,\frac{1}{q}) \nonumber \\
&&=\frac{1}{[2]_q^n} \sum_{k=0}^{n-1}  {n \brack k+1}_q i^{k+1} q^{\frac{k(k+1)}{2}} H_{n-k-1}(x;q)[2]_{\frac{1}{q}}  [k+1]_{\frac{1}{q}} H_{k}(-iqct,\frac{1}{q})(-iqc) \nonumber \\
&&-c  \frac{1}{[2]_q^n} \sum_{k=0}^{n-1}  {n \brack k}_q i^k q^{\frac{k(k-1)}{2}} [2]_q [n-k]_q H_{n-k-1}(x;q) H_k(-iqct,\frac{1}{q}) \nonumber\\
&&= c i^k  \frac{1}{[2]_q^n}\sum_{k=0}^{n-1}  \left( {n \brack k+1}_q i q^{\frac{k(k+1)}{2}} [2]_{\frac{1}{q}}  [k+1]_{\frac{1}{q}} q  - {n \brack k}_q  q^{\frac{k(k-1)}{2}} [2]_q [n-k]_q\right)H_{n-k-1}(x;q) H_k(-iqct,\frac{1}{q})\nonumber\\
&&= 0.\nonumber
\end{eqnarray}
The expression in parenthesis is zero due to $q$-combinatorial formula and $[n]_{\frac{1}{q}}=\frac{[n]_q}{q^{n-1}}$.

\section*{Acknowledgments}
This work was supported by Izmir Institute of
Technology. One of the authors (S. Nalci) was partially supported by TUBITAK scholarship for graduate students.

\end{document}